\documentclass{amsart}
\usepackage{mathrsfs,amssymb,amsfonts,amsmath}
\newtheorem{pro}{Proposition}[section]

\newtheorem{nota}{Remark}[section]

\newcommand{\re}{\mathbb{R}}

\title{Ultraspherical type generating functions for orthogonal polynomials}
\begin{document}
\maketitle
\begin{abstract} 
We characterize, under some technical assumptions and up to a conjecture, probability distributions of finite all order moments with ultraspherical type generating functions for orthogonal polynomials. Our method is based on differential equations and the obtained measures are particular Beta distributions. We actually recover the free Meixner family of probability distributions so that our method gives a new approach to the characterization of free Meixner distributions.  
\end{abstract}
\section{Motivation: Meixner families}
There is a one to one correspondance between probability distributions on the real line and polynomials of a one variable satisfying a three-terms recurrence relation subject to some positivity conditions (\cite{Ism}). That is why in most of the cases, if not all, one tries to characterize probability distributions using generating functions for orthogonal polynomials. Among the famous generating functions are the ones of \emph{exponential type}, that is  if $\mu$ is a probability distribution with a finite exponential moment in a neighborhood of zero 
\begin{equation*}
\int_{\re} e^{zx} \mu(dx) \quad < \quad \infty,
 \end{equation*}  
then 
\begin{equation}\label{GF1}
\psi(z,x) := \sum_{n \geq 0}P_n(x)z^n = \frac{e^{xH(z)}}{\mathbb{E}(e^{XH(z)})}, 
\end{equation}   
where $H$ is analytic around $z = 0$ such that $H(0) = 0, H'(0)=1$, $X$ is a random variable in some probability space $(\Omega, \mathscr{F}, \mathbb{P})$ with law 
$\mu = \mathbb{P} \circ X^{-1}$ and $(P_n)_{n \geq 0}$ is the set of orthogonal polynomials with respect to $\mu$. 
Up to translations and dilations, there are six probability distributions which form the so-called Meixner family referring to its first appearance with J. Meixner (\cite{Meix}). It consists of Gaussian, Poisson, Gamma, negative binomial, Meixner and binomial distributions. This family appeared many times under differents guises (\cite{She}, \cite{Laha}, \cite{AlS}, \cite{Mor}, \cite{Kubo}).\\
Another well known example was first suggested and studied in \cite{Ans} and is given by a \emph{Cauchy-Stieltjes type} kernel. Namely, if $\mu$ is a probability distribution of finite
all order moments, then 
\begin{equation}
\label{GF2}
\psi(z,x) := \sum_{n \geq 0}P_n(x)z^n = \frac{1}{u(z)[f(z) - x]}
\end{equation}
where $u$ and $z \mapsto zf(z)$ are analytic functions around zero such that 
\begin{equation*}
\lim_{z \rightarrow 0} \frac{u(z)}{z} = \lim_{z\rightarrow 0} zf(z) = 1.
\end{equation*}
This family, known as the free Meixner family due to its intimate relation to free probability theory, covers six compactly-supported probability measures too. We refer the reader to 
\cite{Boz}, \cite{Boz1}, \cite{Bry1}, \cite{Kubo2} for more characterizations and more interpretations. The natural $q$-deformation that interpolates the forementioned families for arbitrary $|q| \leq 1$ was defined and studied in \cite{Ans1} and is up to affine transformations the so-called Al-Salam and Chihara family of orthogonal polynomials (\cite{AlS}). Their generating functions is given by an infinite product and is somehow similar to the $q$-exponential function. Another characterization of the last family was recently  given in \cite{Bry}. \\
After this sketchy overview, we suggest another type of generating functions which may be viewed as a generalization of the free Meixner family. It is inspired from the case of Gegenbauer or ultraspherical polynomials for which (\cite{Ism})  
\begin{equation}
\label{UGF}
\sum_{n \geq 0}2^n\frac{(\lambda)_n}{n!}C_n^{\lambda}(x)z^n = \frac{1}{(1- 2zx + z^2)^{\lambda}},\, |x| \leq 1, \lambda > 0,
\end{equation}
where $(\lambda)_n = (\lambda+n-1)\cdots(\lambda+1)\lambda$ and for complex $z$ such that the RHS makes sense and the series in the LHS converges. 
We adopted here the monic normalization for $(C_n^{\lambda})_n$ and henceforth all the polynomials are monic so that they satisfy the normalized recurrence relation
\begin{equation}\label{NRR}
xP_n(x) = P_{n+1}(x) + \alpha_nP_n(x) + \omega_nP_{n-1}(x), \, n \geq 0, P_{-1} := 0, \omega_0 = 1.
\end{equation}
The sequences $(\alpha_n)_{n \geq 0}, (\omega_n)_{n \geq 0}$ are known as the Jacobi-Szeg\"o parameters and $\omega_n > 0$ for all $n$ unless $\mu$ is has a finite support (\cite{Ism}). Moreover, we shall always use these notations for the different families of orthogonal polynomials we shall cross through this paper. \\
It is then natural to adress the problem of characterizing probability measures of finite all order moments, say $\mu_{\lambda}$, such that
\begin{equation}\label{UTGF}
\psi_{\lambda}(z,x) := \sum_{n \geq 0} \frac{(\lambda)_n}{n!}P_n^{\lambda}(x) z^n = \frac{1}{u_{\lambda}(z)(f_{\lambda}(z) - x)^{\lambda}}, \quad \lambda > 0,
\end{equation}  
valid for $x \in \textrm{supp}(\mu_{\lambda})$ and $z$ belongs to a complex open region $S$ near $z=0$ cut from $z=0$ along the negative real axis where $u_{\lambda}, f_{\lambda}$ are analytic with  
\begin{equation*}
(\star) \quad \lim_{z \rightarrow 0} zf_{\lambda}(z) = 1, \quad \lim_{\substack{z \rightarrow 0\\ z \in S}}\frac{u_{\lambda}(z)}{z^{\lambda}} = 1, \,\Im(f(z)) \neq 0,\, z \in S.
\end{equation*} 
By the last assumption, $(f(z) - x)^{\lambda}$ is well defined for all $x \in \textrm{supp}(\mu), z \in S$ and $\lambda > 0$ (the principal determination of the Logarithm is adopted). Moreover, the above limiting conditions imply that $\psi_{\lambda}(z,x)$ tends to $1$ as $z$ tends to $0$ in $S$ for all $x \in \textrm{supp}(\mu_{\lambda})$. We shall say that $\psi_{\lambda}$ is a generating function for orthogonal polynomials of ultraspherical-type referring to ultraspherical polynomials. Without loss of generality, we may assume that $\mu_{\lambda}$ is standard, that is, has a zero mean and a unit variance. Equivalently, if $(\alpha_n^{\lambda})_{n \geq 0}, (\omega_n^{\lambda})_{n \geq 0}$ denote the Jacobi-Szeg\"o parameters of $\mu_{\lambda}$, then one has $\alpha_0^{\lambda} = 0, \omega_1^{\lambda} = 1$. Our strategy is based on the following general claim that was stated without proof in \cite{Boz2} and proved below for the reader's convenience: \\
{\bf Claim}: {\it to a given generating function for orthogonal polynomials $(z,x) \mapsto \psi(z,x)$ associated with a (standard) probability measure $\mu$ satisfying some integrability conditions 
(to be precise later), the measures $\{\mathbb{P}_z\}$ defined by 
\begin{equation*}
\mathbb{P}_z(dx) := \psi(z,x)\mu(dx)
\end{equation*} 
are probability measures such that the mean and the variance of $\mathbb{P}_z$ are polynomials in $z$ of degree $1,2$ respectively.} $\{\mathbb{P}_z\}$ is then referred to as the $\psi$-family of $\mu$ with an at most quadratic variance, referring to both the exponential and the Cauchy-Stieltjes families (\cite{Bry1},\cite{Mor}). 
When $\psi$ is handable enough so that one can perform computations of the first and of the second moments of $\mathbb{P}_z$ independently from the infinite series, one recovers two equations that may be used to solve the problem of characterization of probability measures whose generating function for orthogonal polynomials is given by $\psi$ (or of $\psi$-type). In the case of the Meixner and the free Meixner families, this was noticed in \cite{Boz2}. In the case in hands, if the assumptions in ($\star$) are valid for $z\in S$ together with the assumption ($\star\star$) (see below), we obtain 

\begin{pro} \label{Prop}\noindent
\begin{enumerate}
\item The function $f_{\lambda}$ satisfies for $z \in S$ 
\begin{align}\label{RIC}
Q_2(z)f'_{\lambda}(z) = f_{\lambda}^2(z) - Q_1(z) f_{\lambda}(z) +  R_1(z)
\end{align}
where $Q_2,R_1$ are polynomials of degree $2$ while $Q_1$ is a polynomial of degree $1$. Moreover the coefficients of these polynomials depend only on $\lambda, 
\alpha_1^{\lambda}, \omega_2^{\lambda}$.
\item The function $u_{\lambda}$ is related to $f_{\lambda}$ by 
\begin{equation*}
\frac{u_{\lambda}'(z)}{u_{\lambda}(z)} = \lambda \frac{1-f'_{\lambda}(z)}{f_{\lambda}(z) - \lambda z}. \qquad \qquad \qquad \qquad \blacksquare
\end{equation*}
\end{enumerate}
\end{pro}     
Once we did, we show that if    
\begin{equation}\label{form}
g_{\lambda}(z) := f_{\lambda}(z) - \frac{Q_1(z)}{2} := \frac{E_{\lambda}(z)}{z}
\end{equation}
where $E_{\lambda}$ is assumed to be a polynomial, then $\textrm{deg}(E_{\lambda}) \leq 2$ and this follows from the fact that $Q_2,Q_1,R_1$ are polynomials (terminating series). Next, we investigate under the last assumption the case of symmetric measures. We show that there exist two families of probability measures corresponding to $(C_n^{\lambda})_n$ for 
$\lambda > 0$ and $(C_n^{\lambda-1})_n$ for $\lambda > 1/2, \lambda \neq 1$. We warn the reader to the fact that, though these two families differ from each other by a parameter's translation, their generating functions given by (\ref{UTGF}) are totally different since $a_n^{\lambda}$ depends on $\lambda$ and is fixed for both families. Under the same assumption, there is only one family of non symmetric probability measures corresponding to shifted monic Jacobi polynomials $P_n^{\lambda-1/2,\lambda-3/2}, P_n^{\lambda-3/2,\lambda-1/2}$ for $\lambda > 1/2,\lambda \neq 1$. The discard of the value $\lambda =1$ is needed for the computations since we need to remove factors like $1-\lambda, 1-\lambda^2$. Thus, one deals with this case separately and recovers the free Meixner family for which $\deg(E_\lambda) \leq 1$ too. \\
{\bf Problems}: we do not know if there exists a solution $f_{\lambda}$ for which $E_{\lambda}$ is an entire infinite series. Note that such a solution does not exist when $\lambda=1$. However, we already know that the free Meixner family covers six families of probability distributions (\cite{Boz}) while there are three families for $\lambda \neq 1$ when $E_{\lambda}$ is a polynomial. Is there any intuitive explanation to this difference between both cases or to the degeneracy of the case $\lambda=1$?   

\section{Validity and Proof of the claim}
Write $\psi$ as   
\begin{equation*}
\psi(z,x) = \sum_{n \geq 0}a_nP_n(x)z^n
\end{equation*}
for some fixed sequence $(a_n)_n, x \in \textrm{supp}(\mu)$ and $z$ in a suitable complex domain $D$ near $z=0$ so that the infinite series converge. The integrability conditions we need for the claim to be valid are the finiteness of all order moments of $\mu$ and 
\begin{equation*} 
\int \sum_{n \geq 0} a_n (x^iP_n(x))z^n \mu(dx) =  \sum_{n \geq 0} a_n \int x^iP_n(x)\mu(dx)z^n, \quad i \in \{0,1,2\},\quad (\star\star)
\end{equation*} 
for $z\in D$. 
In fact, for $i=0$, the orthogonality of $P_n$ shows that $\mathbb{P}_z$ is a probability measure for all $z \in D$ (remember that $P_0 = 1$) and together with 
$\alpha_0=0,\omega_1 = 1$ imply  
\begin{eqnarray*}
\sum_{n \geq 0}a_n\int P_{n+1}(x)\mu(dx) z^n &=& 0, \, n \geq 0,  \sum_{n \geq 0}a_n\alpha_n \int P_{n}(x)\mu(dx) z^n = a_0\alpha_0 = 0, \\
\sum_{n \geq 0}a_n\omega_n \int P_{n-1}(x)\mu(dx) z^n &=& a_1\omega_1 = a_1z. 
\end{eqnarray*}
Thus, one gets for $i=1$ after using (\ref{NRR}) 
\begin{equation}\label{M1}
\int \sum_{n \geq 0} a_n (xP_n(x))z^n \mu(dx) = a_1z = \int x\mathbb{P}_z(dx).
\end{equation}
For $i=2$, one uses twice (\ref{NRR}) to get 
 \begin{equation}\label{M2}
\int \sum_{n \geq 0} a_n (x^2P_n(x))z^n \mu(dx) = a_2\omega_2z^2 + a_1\alpha_1 z +1 = \int x^2\mathbb{P}_z(dx)
\end{equation}
and the claim is proved. 

\begin{nota}\label{RE1}
In the case in hands, if $\mu_{\lambda}$ is compactly supported, then the Jacobi-Szeg\"o parameters are bounded thereby one can exchange the infinite sum and integral signs. 
Indeed, by Cauchy-Schwarz inequality
\begin{equation*}
\sum_{n \geq 0} \frac{(\lambda)_n}{n!} \int |x^iP_n(x)|\mu_{\lambda}(dx)|z|^n \leq \left(\int |x|^{2i}\mu(dx)\right)^{1/2} \sum_{n \geq 0}\frac{(\lambda)_n}{n!} ||P_n||\, |z|^n,
\end{equation*}
for $i \in \{0,1,2\}$. Moreover, $||P_n||^2 = \omega_0\dots \omega_{n-1} < c^n$ for some $c>0$ so that Fubini's Theorem applies for $|z| < 1/\sqrt{c}$. As the reader can see, the exchange of the order of integration depends on the sequence $(a_n)_n$ and the growth conditions satisfied by $\mu$. As a matter of fact, if $(a_n)_n$ is fixed, they solely depend on $\mu$ (or in $||P_n||$).
\end{nota}

\section{Proof of Proposition \ref{Prop}}
\subsection{First and second moments}
On the one hand, the integration of both sides of (\ref{UTGF}) with respect to $\mu_{\lambda}$ gives 
\begin{equation*}
u_{\lambda}(z) = \int_{\re}\frac{1}{(f_{\lambda}(z) - x)^{\lambda}}\mu_{\lambda}(dx).
\end{equation*} 
 On the other hand, one gets from (\ref{M1}), (\ref{M2}) and $a_n = (\lambda)_n/n!$
\begin{eqnarray*}
m_1^{\lambda}(z) &:=&  \int x\psi_{\lambda}(z,x)\mu(dx) = \lambda z,\\
m_2^{\lambda}(z) &:=&  \int x^2\psi_{\lambda}(z,x)\mu(dx) = \frac{\lambda(\lambda +1)}{2}\omega_{2}^{\lambda}z^2 +\lambda\alpha_{1}^{\lambda}z + 1.
\end{eqnarray*}
Then, using the elematary operation $x= (x-f(z))+f(z)$, it follows that
\begin{equation*}
m_1^{\lambda}(z) = f(z) - \frac{u_{\lambda,1}(z)}{u_{\lambda}(z)}, \quad u_{\lambda,1}(z) := \int_{\re}\frac{1}{(f(z) - x)^{\lambda-1}}\mu_{\lambda}(dx).
 \end{equation*}
Differentiating with respect to $z \in S$ under the integral sign \footnote{This is justified by the analyticity of $f_{\lambda}$ in $S$ and general properties of generalized Cauchy-Stieltjes transforms, see \cite{Sum} and references therein.} defining $u_{\lambda,1}$, one gets $(1-\lambda) f'(z)u_{\lambda}(z) = (u_{\lambda,1})'(z)$.  
Thus the RHS of  $m_1^{\lambda}(z)$ transforms to:
\begin{equation}\label{E2}
\frac{u_{\lambda}'(z)}{u_{\lambda}(z)} = \lambda \frac{1-f'_{\lambda}(z)}{f_{\lambda}(z) - \lambda z}
\end{equation}
which can be written as
\begin{equation}\label{E2bis}
(u_{\lambda}(z)[f_{\lambda}(z) - \lambda z)])' = (1-\lambda) u_{\lambda}(z)f'_{\lambda}(z).
\end{equation}
For the second moment, use $x^2 = x(x-f(z)) + xf(z)$ to get
\begin{equation}\label{E3}
m_2^{\lambda}(z) = \lambda z f_{\lambda}(z) - \frac{1}{u_{\lambda}(z)} \int_{\re}\frac{x}{(f_{\lambda}(z) - x)^{\lambda-1}}\mu_{\lambda}(dx).
\end{equation} 
Using
\begin{equation*}
\left(\int_{\re}\frac{x}{(f_{\lambda}(z) - x)^{\lambda-1}}\mu_{\lambda}(dx)\right)' = (1-\lambda)f_{\lambda}'(z)\int_{\re}\frac{x}{(f_{\lambda}(z) - x)^{\lambda}}\mu_{\lambda}(dx)
 = \lambda(1-\lambda) zu_{\lambda}(z)f_{\lambda}'(z)
 \end{equation*} 
 $(\ref{E3})$ is rewritten as 
\begin{equation}\label{E3bis}
\left([\lambda z f_{\lambda}(z) - m_2^{\lambda}(z)]u_{\lambda}(z)\right)' = \lambda(1-\lambda)zu_{\lambda}(z) f'_{\lambda}(z). 
\end{equation}

\subsection{A non linear differential equation} 
By the virtue of (\ref{E2bis}), (\ref{E3bis}) implies that
\begin{align*}
\left([\lambda z f_{\lambda}(z) - m_2^{\lambda}(z)]u_{\lambda}(z)\right)' = \lambda z (u_{\lambda}(z)[f_{\lambda}(z) - \lambda z)])'
\end{align*} 
which gives 
\begin{align*}
&[\lambda z f_{\lambda}(z) - m_2^{\lambda}(z)]u'_{\lambda}(z) + [\lambda  f_{\lambda}(z) + \lambda zf'_{\lambda}(z)- (m_2^{\lambda})'(z)]u_{\lambda}(z) 
 \\& = \lambda z [f_{\lambda}(z) - \lambda z]u'_{\lambda}(z) + \lambda z [f'_{\lambda}(z) - \lambda)]u_{\lambda}(z),
 \end{align*}
therefore 
\begin{equation*}
[\lambda^2 z^2 - m_2^{\lambda}(z)] u'_{\lambda}(z) = [(m_2^{\lambda})'(z) - \lambda f_{\lambda}(z) - \lambda^2 z]u_{\lambda}(z).
\end{equation*}
If $\lambda z - m_2^{\lambda}(z) \neq 0$, one gets after the comparison of the last equality to (\ref{E2}) 
\begin{equation*}
\frac{(m_2^{\lambda})'(z) - \lambda f_{\lambda}(z) - \lambda^2 z}{\lambda^2 z^2 - m_2^{\lambda}(z)} = \lambda \frac{1-f'_{\lambda}(z)}{f_{\lambda}(z) - \lambda z}
\end{equation*} 
which shows after elemantary computations that $f_{\lambda}$ satisfies the following non linear first order differential equation:
\begin{align}\label{DE1}
Q_2(z)f'_{\lambda}(z) = f_{\lambda}^2(z) - Q_1(z) f_{\lambda}(z) +  R_1(z)
\end{align}
where 
\begin{eqnarray*}
Q_2(z) &=& \lambda\left[\lambda - \frac{\lambda+1}{2}\omega_2^{\lambda}\right] z^2 - \lambda \alpha_1^{\lambda}z -1,\\ 
Q_1(z) &=& (\lambda+1)\omega_2^{\lambda}z + \alpha_1^{\lambda},\\ 
R_1(z) &=& \frac{\lambda(\lambda+1)}{2}\omega_2^{\lambda} z^2 -1.
\end{eqnarray*}
Setting $g_{\lambda}(z) := f_{\lambda}(z) - [Q_1(z)/2]$, (\ref{DE1}) transforms to 
\begin{equation}
\label{DE2}
Q_2(z)g'_{\lambda}(z) = g_{\lambda}^2(z)  + \tilde{Q}_2(z)
\end{equation} 
where 
\begin{align*}
\tilde{Q}_2(z) &= R_1(z) - \frac{1}{4}[Q_1(z)]^2 - \frac{\lambda+1}{2}\omega_2^{\lambda}Q_2(z) \\
& = [(\lambda+1)\omega_2^{\lambda} -2\lambda]\frac{\lambda^2-1}{4}\omega_2^{\lambda}z^2 + \frac{\lambda^2-1}{2}\alpha_1^{\lambda}\omega_2^{\lambda} z
+ \frac{(\lambda+1)\omega_2^{\lambda}}{2}-1 - \frac{(\alpha_1^{\lambda})^2}{4}.
\end{align*}
Finally, once $g_{\lambda}$ is given, one deduces $f_{\lambda}$ by adding $Q_1/2$ then use $(\ref{E2})$ to derive $u_{\lambda}$.

\section{Some solutions of (\ref{RIC})} From now on, we shall look for solutions of (\ref{RIC}) of the form 
\begin{equation*}
g_{\lambda}(z) := \frac{E_{\lambda}(z)}{z},\quad E_{\lambda}(0) = 1
\end{equation*}
for a second degree polynomial $E_{\lambda}$. In fact, since $z\mapsto zg_{\lambda}(z)$ is analytic around zero, one may always assume that $g_{\lambda}(z)$ has the above form for an entire function $E_{\lambda}$. But if $E_{\lambda}$ is a polynomial of degree $\geq 3$, then all the terms of degree $\geq 3$ will vanish only by equating both sides of (\ref{DE2}). For instance, let 
\begin{equation*}
E_{\lambda}(z)= a_0z^3+ a_1z^2 + a_2z + a_3
\end{equation*} 
and write (\ref{DE2}) as
\begin{equation}
Q_2(z)[zE_{\lambda}'(z) - E_{\lambda}(z)] - E_{\lambda}^2(z)= z^2 \tilde{Q}_2(z).
\end{equation} 
Then by equating terms of degree $6$ is this equation, one easily gets $a_0= 0$ so that $E_{\lambda}$ has degree $2$. For $E_{\lambda}$ a polynomial of degree $4$, start with equating terms of degree $8$ and so on. However, this way of thinking fails or rather become cumbersome when $E_{\lambda}$ is an entire function and the existence of such a solution is open.   

\subsection{A new approach to the Free Meixner family}
Recall that the free Meixner family corresponds to $\lambda=1$ and that it covers six compactly-supported probability distributions given by their Jacobi-Szeg\"o parameters 
(\cite{Boz})
\begin{equation*}
\alpha_n^1 = a, a \in \mathbb{R}, n \geq 1, \quad \omega_n^1 = (1+b), b \geq -1, n \geq 2,
\end{equation*}
where we used the fact that $\mu_1$ has a mean zero $(\alpha_0^1 = 0)$ and a unit variance $(\omega_1^1 = 1)$. Moreover, one has (\cite{Boz1})
\begin{equation*}
f_1(z) = \frac{1+az + (1+b)z^2}{z} \quad \Rightarrow \quad g_1(z) = \frac{(a/2)z +1}{z} = \frac{a}{2} + \frac{1}{z}.
\end{equation*}
But, $\tilde{Q}_2$ reduces to a constant for $\lambda=1$ so that $(\ref{DE2})$ transforms to
\begin{equation*}
[(1-\omega_2^{\lambda})z^2 - \alpha_1^{\lambda}z -1]g'_{\lambda}(z) = g_{\lambda}^2(z) +  (\omega_2^{\lambda}-1) - (\alpha_1^{\lambda})^2/4.
\end{equation*}
It is then an easy exercice to check that $g_1$ satisfies $(\ref{DE2})$ which reads in this case
\begin{equation}
-[bz^2+ az + 1]g'_1(z) = g_1^2(z) + b - a^2/4.
\end{equation} 
We can even prove that $g_1$ as written above is the unique solution of the last differential equation subject to the condition $zg_1(z) \rightarrow 1$ when $z \rightarrow 0$. In fact, writing $g_1(z) = h_1(z) + 1/z$ is some punctured neighborhood of zero where $h_1$ is analytic around zero, simple manipulations show that $h_1$ satisfies
\begin{equation*}
-[bz^2+ az + 1]h'_1(z) = h_1^2(z) - \frac{a^2}{4} + \frac{2}{z}\left(h_1(z) - \frac{a}{2}\right).
\end{equation*}
Taking the limit as $z \rightarrow 0$, one has from the singularity at $z=0$ in the RHS that $h_1(0) = a/2$. Thus, one gets 
\begin{equation*}
-[bz^2+az+1]\sum_{n \geq 1}nc_nz^{n-1} =  \sum_{n \geq 1}c_nz^{n}\left[\sum_{n \geq 1}c_nz^{n} +a\right] + 2\sum_{n \geq 1}c_nz^{n-1}
\end{equation*}
for some sequence $(c_n)_{n \geq 1}$, which makes sense for $z=0$ therefore $c_1 = 0$. Removing $z$ from both sides of the obtained equation then setting $z=0$ will give $c_2=0$, removing $z^2$ and taking $z=0$ gives $c_3 = 0$ and so on. As a result, $h_1(z) = a/2$ and our method gives a new (geometrical) approach to the characterization of free Meixner distributions. 
\begin{nota}
When $\lambda\neq 1$, auxiliary terms show up and $h_{\lambda}$ satisfies 
\begin{equation*}
Q_2(z)h_{\lambda}'(z) = h_{\lambda}^2(z) - \frac{(\alpha_1^{\lambda})^2}{4} +  \frac{2}{z}\left(h_{\lambda}(z) - \lambda\frac{\alpha_1^{\lambda}}{2}\right) + 
\frac{\lambda^2-1}{2}(2-\omega_2^{\lambda})
\end{equation*}
which shows that $h_{\lambda}(0) = \lambda\alpha_1^{\lambda}/2$ while $h_{\lambda}'(0)\neq 0, h_{\lambda}''(0) \neq 0$ in general.
\end{nota}

\section{Symmetric measures: ultraspherical polynomials}
In the sequel, we shall focus on the case $\alpha_{n}^{\lambda} = 0$ for all $n$. This is equivalent to the fact that $\mu_{\lambda}$ is symmetric, that is the image of $\mu_{\lambda}$ by the map $x \mapsto -x$ is still $\mu_{\lambda}$. In this case, one gets by taking $\alpha_1^{\lambda} = 0$  
\begin{eqnarray*}
Q_2(z) &=& \frac{\lambda}{2}[2\lambda - (\lambda+1)\omega_2^{\lambda}] z^2 -1,\\
\tilde{Q}_2(z) &=&  [(\lambda+1)\omega_2^{\lambda} -2\lambda]\frac{\lambda^2-1}{4}\omega_2^{\lambda}z^2  + \frac{(\lambda+1)\omega_2^{\lambda}}{2}-1.
\end{eqnarray*} 
Writing $E_{\lambda}(z) = a_0z^2 + a_1z + a_2$ and equating both sides in (\ref{DE2}), one gets: 
\begin{eqnarray*}
a_2 &=& 1, \\
a_1 &=& 0,\\
-3a_0 - \frac{\lambda}{2}[2\lambda - (\lambda+1)\omega_2^{\lambda}] &=&  \frac{(\lambda+1)\omega_2^{\lambda}}{2}-1,\\ 
-a_0^2 + a_0 \frac{\lambda}{2}[2\lambda - (\lambda+1)\omega_2^{\lambda}]  &=& [(\lambda+1)\omega_2^{\lambda} -2\lambda]\frac{\lambda^2-1}{4}\omega_2^{\lambda}.
\end{eqnarray*}
The third equation gives 
\begin{equation*}
a_0 = \frac{(1-\lambda^2)(2-\omega_2^{\lambda})}{6}.
\end{equation*}
Hence, it remains to check when the above $a_0$ satisfies the fourth equation. Since the case $\lambda=1$ is known, we assume $\lambda \neq 1$ so that one removes the term 
$(1-\lambda^2)$ in the above equalities.  Substituting $a_0$ in the fourth equation, one sees that $\omega_2^{\lambda}$ satisfies
\begin{equation*}
-(\lambda+1)(\lambda +2) (\omega_2^{\lambda})^2 + (4\lambda^2+ 6\lambda - 1)\omega_2^{\lambda} + (1- 4\lambda^2) = 0.
\end{equation*}
What is quite interesting and even surprising, that though this polynomial looks complicated, its descriminant is equal $9$ so that there are two solutions given by 
\begin{equation*}
 \omega_{2,1}^{\lambda} = \frac{2\lambda+1}{\lambda+2} \qquad \omega_{2,2}^{\lambda} = \frac{2\lambda-1}{\lambda+1}
 \end{equation*} 
where for the second value, we consider $\lambda > 1/2$ in order to avoid finitely-supported probability measures and signed measures. 
As a result, 
\begin{equation*}
a_0 = \frac{1-\lambda^2}{2(\lambda +2)} , \qquad a_0 = \frac{1-\lambda^2}{2(\lambda +1)} = \frac{1-\lambda}{2}.
\end{equation*}
Thus 
\begin{equation*}
f_{\lambda}(z) = \frac{1+\lambda}{2}z +\frac{1}{z},\quad f_{\lambda}(z) = \frac{\lambda}{2}z + \frac{1}{z},  
 \end{equation*}
 and from (\ref{E2}) 
 \begin{equation*}
\frac{u_{\lambda}'(z)}{u_{\lambda}(z)} = \frac{\lambda}{z},\quad \frac{u_{\lambda}'(z)}{u_{\lambda}(z)} = \lambda \frac{z^2+1-(\lambda/2) z^2}{z(1-(\lambda/2)z^2)}.
 \end{equation*}
Finally
 \begin{equation*}
u_{\lambda}(z) = z^{\lambda},\, \lambda > 0, \lambda \neq 1,\quad u_{\lambda}(z)= \frac{z^{\lambda}}{1-(\lambda/2)z^2}, \lambda > 1/2, \lambda \neq 1,
 \end{equation*}
for $z \in S$. Note that $S$ is easily described: in fact $f_{\lambda}$ is not real outside the real line and the circle $|z| < 2/(1+\lambda)$ or $|z| < 2/\lambda$ respectively. Moreover the $\mu_{\lambda}$ is compactly-supported as we shall see below, so that $(\star\star)$ is satisfied in a ball centered at the origin (see remark \ref{RE1}).    

\subsection{Ultraspherical polynomials: symmetric Beta distributions}
The value $\omega_{2,1}^{\lambda}$ corresponds to the ultraspherical polynomials. However, in order to fit into our setting, one has to consider the monic Gegenbauer polynomials,
say $ \tilde{C_n^{\lambda}}$, which are orthogonal with respect to the standard Beta distribution 
\begin{equation*}
c_{\lambda}(1-x^2/[2(1+\lambda)])^{\lambda -1/2}dx, \quad x \in [\pm \sqrt{2(1+\lambda)}]
\end{equation*}
for some normalizing constant $c_{\lambda}$. They are given by 
\begin{equation*}
\tilde{C_n^{\lambda}}(x) = (\sqrt{2(1+\lambda)})^n C_n^{\lambda}\left(\frac{x}{\sqrt{2(1+\lambda)}}\right). 
\end{equation*}
Now, it is easy to see from (\ref{UGF}) that 
\begin{align*}
\sum_{n \geq 0}\frac{(\lambda)_n}{n!} \tilde{C_n^{\lambda}}(x) z^n &= 
\sum_{n \geq 0}2^n\frac{(\lambda)_n}{n!} C_n^{\lambda}\left(\frac{x}{\sqrt{2(1+\lambda)}}\right) \left(\frac{\sqrt{1+\lambda}z}{\sqrt 2}\right)^n 
\\& = \frac{1}{(1- zx + (1+\lambda)z^2/2)^{\lambda}} 
\\&=  z^{-\lambda}\left[\frac{1+ (1+\lambda)z^2/2}{z} - x\right]^{-\lambda} = \frac{1}{u_{\lambda}(z)(f_{\lambda}(z) - x)^{\lambda}}. \qquad \qquad \blacksquare
\end{align*}

For $\omega_{2,2}^{\lambda}$,  $\psi_{\lambda}$ is written as:
\begin{equation*}
\psi_{\lambda}(z,x)  = \frac{1-(\lambda/2) z^2}{z^{\lambda}(\lambda z/2 + 1/z -x)^{\lambda}} = \frac{1-(\lambda/2) z^2}{(\lambda z^2/2 + 1 - zx)^{\lambda}}
\end{equation*}
and we claim that $P_n^{\lambda} = \tilde{C_n}^{\lambda-1}$ for all $n$ and all $\lambda> 1/2,\lambda \neq 1$. In fact,  
\begin{align*}
\sum_{n \geq 0}\frac{(\lambda)_n}{n!} \tilde{C_n}^{\lambda-1}(x) z^n &= \sum_{n \geq 0}\frac{\lambda +n-1}{\lambda-1}\frac{(\lambda-1)_n}{n!} \tilde{C_n}^{\lambda-1}(x) z^n
\\& = \frac{1}{(\lambda-1)z^{\lambda-2}}\partial_z \sum_{n \geq 0}\frac{(\lambda-1)_n}{n!} \tilde{C_n}^{\lambda-1}(x) z^{n+\lambda-1}
\\& = \frac{1}{(\lambda-1)z^{\lambda-2}} \partial_z \left[\frac{z}{1-zx +\lambda z^2/2}\right]^{\lambda-1}
\\& =  \frac{1-(\lambda/2) z^2}{( 1 - zx + \lambda z^2/2)^{\lambda}}
\end{align*}  
as the reader may easily check.  $\hfill \blacksquare$
\section{non-symmetric probability measures: Jacobi polynomials}
Henceforth, we suppose that $\alpha_1^{\lambda} \neq 0,\lambda \neq 1$ and we will show that there is only one family of probability measures subject to 
\begin{equation*}
g_{\lambda}(z) = \frac{a_0z^2+a_1z + a_2}{z}.
\end{equation*}
Then, we get the following equations
\begin{eqnarray*}
a_2  &=& 1, \\
a_1 &=& \frac{\lambda \alpha_1^{\lambda}}{2} \neq 0,\\
-3a_0 - \frac{\lambda}{2}[2\lambda - (\lambda+1)\omega_2^{\lambda}] - a_1^2 &=&  \frac{(\lambda+1)\omega_2^{\lambda}}{2}-1 - \frac{(\alpha_1^{\lambda})^2}{4},\\ 
-a_0\alpha_1\lambda -2a_0a_1 &=&  \frac{\lambda^2-1}{2}\alpha_1^{\lambda}\omega_2^{\lambda}, \\
-a_0^2 + a_0 \frac{\lambda}{2}[2\lambda - (\lambda+1)\omega_2^{\lambda}]  &=& [(\lambda+1)\omega_2^{\lambda} -2\lambda]\frac{\lambda^2-1}{4}\omega_2^{\lambda}.
\end{eqnarray*}
From the second, third and fourth equations, it follows that
\begin{equation*}
a_0 = \frac{1-\lambda^2}{6}\left[\frac{(\alpha_1^{\lambda})^2}{2} + 2- \omega_2^{\lambda}\right] = \frac{1-\lambda^2}{4\lambda}\omega_2^{\lambda}.
\end{equation*}
Actually, this gives a constraint on $\lambda,\alpha_1^{\lambda}, \omega_2^{\lambda}$:
\begin{equation}
\label{Cons}
\left(\frac{(\alpha_1^{\lambda})^2}{2} + 2\right)\lambda = \left(\lambda + \frac{3}{2}\right)\omega_2^{\lambda}.
\end{equation}
Substituting $a_0$ by $(1-\lambda^2)\omega_2^{\lambda}/(4\lambda)$ and removing $(1-\lambda^2)$, the fifth equation becomes
\begin{equation*}
-\frac{1-\lambda^2}{16\lambda^2}(\omega_2^{\lambda})^2 + \frac{\omega_2^{\lambda}}{8}[2\lambda - (\lambda+1)\omega_2^{\lambda}] 
= [2\lambda - (\lambda+1)\omega_2^{\lambda}] \frac{\omega_2^{\lambda}}{4}.
\end{equation*}
In the non degenerate case $\omega_2^{\lambda} \neq 0$, 
\begin{equation*}
\omega_2^{\lambda} = \frac{4\lambda^3}{2\lambda^3 + 3\lambda^2 -1}.
\end{equation*}
But $-1$ is a double root of the polynomial in the denominator so that
\begin{equation*}
\omega_2^{\lambda} = \frac{2\lambda^3}{(\lambda + 1)^2(\lambda -1/2)},
\end{equation*}
which is positive for $\lambda > 1/2$. Finally, one deduces from (\ref{Cons}) that 
\begin{eqnarray*}
(\alpha_1^{\lambda})^2 &=&  2\left[\frac{(2\lambda +3)\lambda^2}{(\lambda + 1)^2(\lambda -1/2)} - 2\right] = \frac{2}{(\lambda + 1)^2(\lambda -1/2)} > 0,\\
a_0 &=& = \frac{(1-\lambda^2)\lambda^2}{2(\lambda + 1)^2(\lambda -1/2)} = \frac{(1-\lambda)\lambda^2}{(\lambda + 1)(2\lambda -1)}. 
\end{eqnarray*}
It follows that 
\begin{align*}
f_{\lambda}(z) &= \frac{a_0z^2+a_1z + a_2}{z} + \frac{(1+\lambda)\omega_2^{\lambda}z + \alpha_1^{\lambda}}{2} 
\\& = \frac{1}{z}\left[\left(\frac{1-\lambda}{2\lambda} + 1\right)\frac{1+\lambda}{2}\omega_2^{\lambda}z + \frac{\lambda+1}{2}\alpha_1^{\lambda} + 1 \right]
\\&= \frac{1}{z}\left[\frac{\lambda^2}{2\lambda-1}z^2 \pm \frac{1}{\sqrt{2\lambda-1}}z + 1\right]
\end{align*} 
and 
\begin{equation*}
\frac{u'_{\lambda}(z)}{u_{\lambda}(z)} = \frac{\lambda}{z}\left[1- \frac{(\lambda-1)^2}{2\lambda-1}z^2\right]
\left[\frac{\lambda(1-\lambda)}{2\lambda -1} z^2   \pm \frac{1}{\sqrt{2\lambda-1}}z +1\right]^{-1}.
\end{equation*}
The descriminant of the polynomial 
\begin{equation*}
\frac{\lambda(1-\lambda)}{2\lambda -1} z^2   \pm \frac{1}{\sqrt{2\lambda-1}}z +1
\end{equation*}
is easily seen to be: 
\begin{equation*}
\frac{1}{2\lambda-1} -\frac{4\lambda(1-\lambda)}{2\lambda-1} = 2\lambda -1 > 0.
\end{equation*}
It follows that, when $\alpha_1^{\lambda} > 0$, the roots are given by 
\begin{eqnarray*}
z_1 =-\frac{\sqrt{2\lambda-1}}{\lambda},\quad z_2 = -\frac{\sqrt{2\lambda-1}}{1-\lambda}.
\end{eqnarray*}
Writing
\begin{equation*}
1- \frac{(\lambda-1)^2}{2\lambda-1}z^2 = -\frac{(\lambda-1)^2}{2\lambda-1}\left[z + \frac{\sqrt{2\lambda-1}}{1-\lambda}\right]\left[z - \frac{\sqrt{2\lambda-1}}{1-\lambda}\right],
\end{equation*}
one gets
\begin{align*}
\frac{u'_{\lambda}(z)}{u_{\lambda}(z)} = \frac{\lambda-1}{z}\left[z + \frac{\sqrt{2\lambda-1}}{\lambda-1}\right]\left[z + \frac{\sqrt{2\lambda-1}}{\lambda}\right]^{-1} 
= \frac{\lambda}{z} - \frac{1}{z + \sqrt{2\lambda -1}/\lambda}.
\end{align*}
As a result
\begin{equation*}
u_{\lambda}(z) =  \frac{\sqrt{2\lambda-1}}{\lambda}\frac{z^{\lambda}}{z+ \sqrt{2\lambda-1}/\lambda}
\end{equation*}
and the generating function is written as 
\begin{equation} \label{GF3}
\psi_{\lambda}(z,x)  = \frac{\lambda}{\sqrt{2\lambda-1}} \left[z+ \frac{\sqrt{2\lambda-1}}{\lambda}\right]
\left[1 -z\left(x- \frac{1}{\sqrt{2\lambda-1}}\right) +\frac{\lambda^2}{2\lambda-1}z^2 \right]^{-\lambda}.
\end{equation}

In the case $\alpha_1^{\lambda} < 0$, similar computations yield
\begin{equation*}
u_{\lambda}(z) =   -\frac{\sqrt{2\lambda-1}}{\lambda}\frac{z^{\lambda}}{z- \sqrt{2\lambda-1}/\lambda}
\end{equation*}
and 
\begin{equation*}
\psi_{\lambda}(z,x)  = -\frac{\lambda}{\sqrt{2\lambda-1}} \left[z- \frac{\sqrt{2\lambda-1}}{\lambda}\right]
\left[1 - z\left(x+ \frac{1}{\sqrt{2\lambda-1}}\right) + \frac{\lambda^2}{2\lambda-1}z^2\right]^{-\lambda}.
\end{equation*}

\subsection{Orthogonality measures: special Jacobi polynomials}
We will show that $P_n^{\lambda}$ is a shifted monic Jacobi polynomial with parameters depending on $\lambda$. To proceed, recall that (\cite{Koe}) the monic Jacobi polynomials $p_n^{\alpha,\beta}$ are orthogonal with respect to the Beta distribution with density function given by 
\begin{equation*}
c_{\alpha,\beta}(1-x)^{\alpha}(1+x)^{\beta} {\bf 1}_{[-1,1]}(x),\quad \alpha,\beta > -1, 
\end{equation*}
for some normalizing constant $c_{\alpha,\beta}$ and that the non monic Jacobi polynomials $P_n^{\alpha,\beta}$ are related to $p_n^{\alpha,\beta}$ as
\begin{equation*}
P_n^{\alpha,\beta}(x) = \frac{(n+\alpha + \beta +1)_n}{2^n n!}p_n^{\alpha,\beta}(x) =  \frac{(\alpha + \beta +1)_{2n}}{(\alpha+\beta+1)_n 2^n n!}p_n^{\alpha,\beta}(x).
\end{equation*}
We will show that 
\begin{equation*}
P_n^{\lambda}(x) = \left[\frac{2\lambda}{\sqrt{2\lambda-1}}\right]^n p_n^{\lambda-1/2,\lambda - 3/2}\left(\frac{\sqrt{2\lambda-1}x -1}{2\lambda}\right).
\end{equation*}
when $\alpha_1^{\lambda} > 0$ and
\begin{equation*}
P_n^{\lambda}(x) = \left[\frac{2\lambda}{\sqrt{2\lambda-1}}\right]^n p_n^{\lambda-3/2,\lambda - 1/2}\left(\frac{\sqrt{2\lambda-1}x +1}{2\lambda}\right).
\end{equation*}
when $\alpha_1^{\lambda} < 0$.  Before proceeding, note that both cases are related using $P_n^{\alpha,\beta}(x) = (-1)^nP_n^{\beta,\alpha}(-x)$ (\cite{Ism}): 
\begin{equation*}
P_n^{\lambda-3/2,\lambda - 1/2}\left(\frac{\sqrt{2\lambda-1}x +1}{2\lambda}\right) = (-1)^nP_n^{\lambda-1/2,\lambda - 3/2}\left(\frac{\sqrt{2\lambda-1}(-x) -1}{2\lambda}\right)
\end{equation*}
so that their generating functions are the same up to the transformation $(z,x) \mapsto (-z,-x)$. Moreover the orthogonality measures are given by 
\begin{eqnarray*}
\mu_{\lambda}(dx) &=& c_{\lambda}\left(1 - \frac{\sqrt{2\lambda-1}x -1}{2\lambda}\right)^{\lambda-1/2}\left(1 + \frac{\sqrt{2\lambda-1}x -1}{2\lambda}\right)^{\lambda-3/2}dx,\\
\mu_{\lambda}(dx) &=& c'_{\lambda}\left(1 - \frac{\sqrt{2\lambda-1}x +1}{2\lambda}\right)^{\lambda-1/2}\left(1 + \frac{\sqrt{2\lambda-1}x +1}{2\lambda}\right)^{\lambda-3/2}dx,
\end{eqnarray*}
for some normalizing constants $c_{\lambda}, c'_{\lambda}$ and for
\begin{eqnarray*}
x &\in&  \left[\frac{1-2\lambda}{\sqrt{2\lambda-1}}, \frac{1+2\lambda}{\sqrt{2\lambda-1}}\right], \\
x &\in&  \left[-\frac{1+2\lambda}{\sqrt{2\lambda-1}}, \frac{2\lambda-1}{\sqrt{2\lambda-1}}\right]
\end{eqnarray*}
respectively. \\
Now, we proceed to the proof of our claim and we consider the case $\alpha_1^{\lambda} > 0$. To this end, we need (\cite{Koe}) 
\begin{align*}
\frac{1}{(1+t)^{\alpha+\beta+1}}{}_2F_1\left(\substack{\displaystyle \frac{\alpha +\beta +1}{2}, \frac{\alpha +\beta +2}{2} \\ \displaystyle \beta+1}; \frac{2(y+1)t}{(1+t)^2}\right)
&= \sum_{n \geq 0} \frac{(\alpha+\beta+1)_n}{(\beta+1)_n}P_n^{\alpha,\beta}(y) t^n
\\& = \sum_{n \geq 0} \frac{(\alpha+\beta+1)_{2n}}{(\beta+1)_n n!}p_n^{\alpha,\beta}(y) \left(\frac{t}{2}\right)^n 
\end{align*} 
for $|t| < 1, |y| < 1$, where ${}_2F_1$ is the Gauss hypergeometric function (\cite{Ism}). Substituting $(\alpha,\beta)$ by $(\lambda-1/2, \lambda-3/2)$, 
then $(\alpha +\beta +1)/2 = \lambda-1/2 = \beta+1$ so that   
\begin{equation*}
{}_2F_1\left(\substack{\displaystyle \frac{\alpha +\beta +1}{2}, \frac{\alpha +\beta +2}{2} \\ \displaystyle \beta+1}; \frac{2(y+1)t}{1+t^2}\right) = 
{}_1F_0\left(\lambda;  \frac{2(y+1)t}{(1+t)^2}\right) =  \left(1- \frac{2(y+1)t}{(1+t)^2}\right)^{-\lambda},   
\end{equation*}
where we used that ${}_1F_0(\lambda, y) = (1-y)^{-\lambda}$ for $|y| < 1$ (\cite{Ism}). Thus 
\begin{equation*}
\frac{1}{(1+t)^{\alpha+\beta+1}}{}_2F_1\left(\substack{\displaystyle \frac{\alpha +\beta +1}{2}, \frac{\alpha +\beta +2}{2} \\ \displaystyle \beta+1}; \frac{2(y+1)t}{(1+t)^2}\right)
= \frac{1+t}{[1+t^2 - 2ty]^{\lambda}}.
\end{equation*}
Now use the Gauss duplication formula (\cite{Ism})
\begin{equation*}
\sqrt{\pi} \Gamma(2a) = 2^{2a-1}\Gamma(a)\Gamma(a+1/2), \quad a > 0,
\end{equation*}
to see that
\begin{equation*}
 \frac{(\alpha+\beta+1)_{2n}}{(\beta+1)_n}  = \frac{(2\lambda -1)_{2n}}{(\lambda-1/2)_n} = 2^{2n} (\lambda)_n.
\end{equation*}
As a result, 
\begin{equation*}
 \sum_{n \geq 0} \frac{(\lambda)_n}{n!}p_n^{\alpha,\beta}(y) (2t)^n =  \frac{1+t}{[1+t^2 - 2ty]^{\lambda}}.
 \end{equation*}
 It finally remains to substitute in the last equality
\begin{equation*}
y = \frac{\sqrt{2\lambda-1}x -1}{2\lambda}, \quad  t = \frac{\lambda}{\sqrt{2\lambda-1}}z
 \end{equation*}
for small $z$ to see that it is nothing but (\ref{GF3}) and the claim follows.  $\hfill \blacksquare$

\end{document}